\documentclass{amsart}

\usepackage {epsfig}
\usepackage {graphicx}

 \newtheorem{thm}{Theorem}
 
 \newtheorem{lemma}[thm]{Lemma}
 
 \newtheorem{theorem}[thm]{Theorem}
 \theoremstyle{definition}
 
 \theoremstyle{remark}

 \newcommand{\Real}{\mathbb{R}}
 \newcommand{\Complex}{\mathbb{C}}
 \newcommand{\Field}{\mathbb{K}}

\usepackage{graphicx}
\begin{document}

\title[Solution to Rubel's question]
{Solution to Rubel's question about differentially algebraic
dependence on initial values}

\author{Guy Katriel}

\address{}

\email{haggaik@wowmail.com}

\subjclass{}

\keywords{algebraic differential equation, differentially
algebraic functions. Mathematics Subject Classification: 34A09
(12H05)}
\date{}

\dedicatory{}

%%% ----------------------------------------------------------------------

\begin{abstract}
We prove that, for generic systems of polynomial differential
equations, the dependence of the solution on the initial
conditions is not differentially algebraic. This answers, in the
negative, a question posed by L.A. Rubel.

\end{abstract}

%%% ----------------------------------------------------------------------
\maketitle
%%% ----------------------------------------------------------------------

\section{introduction}

This work answers a question posed by L.A. Rubel in \cite{rubel},
a paper which presents many research problems on differentially
algebraic functions (see also \cite{rubelf}). We recall that an
analytic function $f(z)$ (defined on some open subset of either
$\Real$ or $\Complex$) is called {\it{Differentially Algebraic}}
(DA) if it satisfies some differential equation of the form
\begin{equation}
\label{da} Q(z,f(z),f'(z),...,f^{(n)}(z))=0,
\end{equation}
for all $z$ in its domain, where Q is a nonzero polynomial of
$n+2$ variables. A function of several variables is DA if it is DA
in each variable separately when fixing the other variables. It is
known that a function is DA precisely precisely when it is
computable by a general-purpose analog computer \cite{pour}.

Problem 31 of \cite{rubel} asks:

\begin{quote}
{\emph{Given a `nice' initial value problem for a system of
algebraic equations in the dependent variables $y_1,...,y_n$, must
$y_1(x_0)$ be differentially algebraic as a function of the
{\bf{initial conditions}}, for each $x_0$?}}
\end{quote}

Rubel adds: ``we won't say more about what `nice' means except
that the problem should have a unique solution for each initial
condition in a suitable open set".

In other words, let
$$(P)\;\;\; y_k'=p_k(y_1,...,y_m),\;\;  k=1...m$$
be a system of differential equations, with $p_k$ polynomials.
Denote by $y_k(r_1,...,r_m;x)$ the solutions of $(P)$ with initial
conditions
\begin{equation}
\label{init} y_k(0)=r_k,\;\; k=1...m.
\end{equation}

The question is whether it is true that the dependence of $y_i$ on
$r_j$ (fixing $x$ and the other initial conditions $r_k$, $k\neq
j$) is DA.

Here we first answer the question in the {\it{negative}} by
constructing an initial value problem:
\begin{eqnarray}
\label{spec0}
y_{1}'=y_{2}-y_{1},\;\; y_{1}(0)=r_1\\
\label{spec1} y_{2}'=y_{2}(y_{2}-y_{1}),\;\; y_{2}(0)=r_2
\end{eqnarray}
and proving that, for most $x\in \Real$, both $y_1(0,r_2;x)$ and
$y_2(0,r_2;x)$ are {\it{not}} DA with respect to $r_2$.

Moreover, we prove that this phenomenon is in fact {\it{generic}}
in the class of nonlinear polynomial differential systems $(P)$
(and thus the answer to Rubel's question is {\it{very}} negative).
To state this precisely, let ${\bf{S}}(m,d)$ denote the set of all
polynomial systems $(P)$ of size $m$, where the polynomials $p_i$
are of degree at most $d$. By identifying the coefficients of the
various monomials occuring in the $p_i$'s with coordinates,
${\bf{S}}(m,d)$ is a finite dimensional vector-space, so it
inherits the standard topology and measure.

In all our statements the term {\it{generic}} will have the
following meaning: when $U\subset \Real^N$ is an open set and
$V\subset U$, we will say that $V$ is {\bf{generic}} in $U$ if its
complement in $U$ is {\it{both of measure $0$ and of first Baire
category}}.

For each system $P\in {\bf{S}}(m,d)$ and initial conditions
$r=(r_1,...,r_m)\in \Real^m$, let $I(P,r)$ denote the largest
interval containing $0$ for which the solutions of the
initial-value problem (P),(\ref{init}) are defined.  Let
\begin{equation}
\label{lambda}
\Lambda(P)=\{(x,r)\;|\; r\in \Real^m,\; x\in
I(P,r)\}
\end{equation}
($\Lambda(P)$ is the domain of definition of the functions
$y_i(r_1,...,r_m;x)$ - we note that standard existence theory of
ODE's implies that it is an open set; it is also easy to see that
it is connected). For each $1\leq j \leq m$ let
$$\Lambda_j(P)=\{(x,r_1,...,r_{j-1},r_{j+1},...,r_m)\;|\; \exists r_j
\;\; such\; that\; (x,r_1,...,r_m)\in \Lambda(P) \}.$$

\begin{theorem}
\label{main} Assume $m\geq 2$ and $d\geq 2$. For generic $P\in
{\bf{S}}(m,d)$, we have:\newline For any $1\leq i,j \leq m$, and
for generic choice of ($x,r_1,...,r_{j-1},r_{j+1},...,r_{m})\in
\Lambda_j(P)$, the function
\begin{equation}
\label{fi} f_{ij}(z)=y_i(r_1,...,r_{j-1},z,r_{j+1},...,r_m;x)
\end{equation}
is {\it{not}} DA.
\end{theorem}

We note that the construction of a specific example
(\ref{spec0})-(\ref{spec1}) is essential for our proof of the
`generic' result of theorem \ref{main}.

\vspace{0.4cm} The significance of the problem posed by Rubel is
brought out when we note that the class of DA functions is a very
``robust" one:

\noindent (i) The class is closed under many of the
constructions of analysis: algebraic operations, composition of
functions, inversion, differentiation and integration.

\noindent
(ii) The components of an analytic solution of a
{\it{system}} of algebraic differential equations are themselves
DA \cite{rubel4}.

\noindent (iii) A solution of a differential equation
$R(z,f(z),f'(z),...,f^{(n)}(z))=0$, where $R$ is a DA function, is
itself DA \cite{rubel2}. We note a special case $n=0$ of this
result, which can be called `the DA implicit-function theorem',
which we shall have occasion to use later:

\begin{lemma}
\label{implicit} Assume $R(w,z)$ is a DA function, and $f(z)$ is a
real-analytic function satisfying
$$R(f(z),z)=0$$
for all $z$ in some interval. Then $f$ is DA on this interval.
\end{lemma}

For the above reasons, most of the transcendental functions which
are encountered in ``daily life'' are DA (notable exceptions are
the Gamma function and Riemann Zeta function, see \cite{rubel3}).
It is then natural to wonder whether functions obtained by the
construction of looking at the dependence of the ``final value''
on the ``initial value'' (in other terms, the components of the
``Poincar\'e map") of an algebraic differential equation are also
DA. Our results show that this is almost always not so.

\vspace{0.4cm} The following remarks show that the result of
theorem \ref{main} is in some sense close to optimal:

\noindent (i) The condition $d\geq 2$ cannot be removed. Indeed if
$d=0$ the solutions of the system are linear functions in $x$, and
if $d=1$ the solutions are linear combinations of exponentials in
$x$, so the dependence on initial conditions is certainly DA.

\noindent (ii) The condition $m\geq 2$ cannot be removed. Indeed
if $m=1$ then we are dealing with an equation of the form
$$y'=p(y),$$
which is solved in terms of elementary functions, which implies
that the dependence on the initial condition is DA.

\noindent (iii) The restriction to generic systems cannot be
removed, since, for example, when, for each $k$, $p_k$ depends
only on $y_k$, we are back to a decoupled system of equations of
the same form as in (ii) above, so that again we have DA
dependence on initial conditions. However, it might still be
possible to obtain a stronger statement as to the size of the set
of `exceptional' systems. For example, we do not know whether the
set of systems $P\in {\bf{S}}(m,d)$ for which all the functions
$f_{ij}$ are DA is {\it{nowhere dense}} in ${\bf{S}}(m,d)$.

\vspace{0.4cm} In section \ref{alt} we prove a general result
which underlies the proof of our `genericity' results. In section
\ref{nonda} we construct some explicit examples of functions which
we prove are not DA. In section \ref{const} we show that the
functions constructed in section \ref{nonda} arise in the solution
of the system (\ref{spec0})-(\ref{spec1}), and this is used to
prove that the solutions of this system are not DA in the initial
conditions. Finally, in section \ref{proof} we use the specific
example constructed in section \ref{const}, together with the
general result of section \ref{alt}, to prove theorem \ref{main}.

\section{The Alternative Lemma}
\label{alt}

In this section we present a result which we term the
``alternative lemma", which is used several times in the arguments
of the following sections. This result says that if we have a
parametrized family of analytic functions (with the dependence on
the parameters also analytic), then either: (I) all the functions
in the family are DA, or (II) generic functions in the family are
{\it{not}} DA. Thus, to show that a {\it{generic}} function in a
certain family is {\it{not}} DA, it is sufficient to prove that
{\it{one}} of the functions in the family is {\it{not}} DA.

The proof of the lemma is based on two facts:

\noindent (i) The fact that an analytic function on a connected
open set cannot vanish on a `large' subset without vanishing
identically.

\noindent (ii) The Gourin-Ritt theorem \cite{ritt} which says that
any DA function $f(z)$ in fact satisfies a differential equation
of the form (\ref{da}) where Q is a polynomial with {\it{integer
coefficients}}.

We state the next lemma in both `real' and `complex' forms, since
we shall have occasions to use both.

\begin{lemma}
\label{alternative} Let $\Field =\Real$ or $\Complex$. And let
$F:\Omega \rightarrow \Real$, where $\Omega\subset \Real^p \times
\Field$ is open and connected, be an analytic function. For each
$u\in \Real^{p}$, define
$$\Omega(u)=\{ v\in \Field \; | \; (u,v)\in \Omega \}$$
And let
$$\Omega'=\{ u\in \Real^p \; | \; \Omega(u)\neq \emptyset \;\}.$$
Then the following alternative holds: either

\noindent (I) For every $u\in \Omega'$, $F(u,.)$ is DA on
$\Omega(u)$. In fact, all the functions $F(u,.)$, $u\in \Omega'$,
satisfy the {\it{same}} differential equation.

\noindent
or:

\noindent (II) For generic $u\in \Omega'$, $F(u,.)$ is {\it{not}}
DA on $\Omega(u)$.
\end{lemma}

\noindent {\sc{proof:}}  Let {\bf{Q}} be the set of all
polynomials with integer coefficients. {\bf{Q}} is a countable
set. For each $Q\in {\bf{Q}}$, let $\Omega'(Q)$ be the set of
$u\in \Omega'$, for which $f=F(u,.)$ satisfies (\ref{da}) on
$\Omega(u)$. We claim that either $\Omega'(Q)$ is of measure 0 or
$\Omega'(Q)=\Omega'$. Assume that $\Omega'(Q)$ has positive
measure in $\Real^p$. We choose open balls $B_1\subset \Omega'$
and $B_2\subset \Field$ such that (i) $B_1 \times B_2\subset
\Omega$ and (ii) $\Omega'(Q)\cap B_1$ has positive measure. Then
the real-analytic function
$$g(u,v)=Q(v, F(u,v), D_v F(u,v),...,D_v^{(n)}F(u,v))$$
vanishes on the set $(\Omega'(Q)\cap B_1)\times B_2$, which is of
positive measure. By analyticity and the fact that $\Omega$ is
connected, this implies that $g$ vanishes throughout $\Omega$,
which means that alternative (I) holds. So if alternative (I) does
not hold, $\Omega'(Q)$ must be of measure 0. Since it is easy to
see that $\Omega'(Q)$ is a relatively closed set in $\Omega'$, the
fact that it has measure $0$ implies that it is nowhere dense.
Since this is true for any $Q\in {\bf{Q}}$, the countable union
$$K=\cup_{Q\in {\bf{Q}}}{\Omega'(Q)},$$
is a set of measure 0 and first category. By the Gourin-Ritt
theorem, for all $u$ outside $K$, $F(u,.)$ is {\it{not}} DA. Thus
(II) holds.

\section{Some functions which are not differentially algebraic}
\label{nonda}

In this section we define some new functions and prove that they
are not DA. These results will be used in our construction of an
explicit differential equation with non-DA dependence on initial
conditions in section \ref{example}.

We define the function $H(c)$ (of $c\in\Complex$) by:
\begin{equation}
\label{defh} H(c)=\int_0^{\infty}{\frac{du}{ce^{u}-u}}
\end{equation}

For $c\in \Complex$ to belong to the domain of definition $D_H$ of
$H$ we need to ensure that the denominator $ce^{u}-u$ does not
vanish for any $u\geq 0$. It is then easy to check that
\begin{equation}
\label{domain}
 D_H=\Complex - [0,\frac{1}{e}].
\end{equation}

\begin{lemma}
\label{nda1} The function $H(c)$ is {\it{not}} DA.
\end{lemma}

\noindent {\sc{proof:}} We define $h(z)=H(\frac{1}{z})$. We shall
show that $h$ is {\it{not}} DA, which implies that $H$ is not DA.
We expand $h$ in a power series (which converges for $|z|<e$):
\begin{eqnarray*}
h(z)=\int_0^{\infty}{\frac{ze^{-u} du}{1-zue^{-u}}}
=\int_0^{\infty}{ze^{-u}\sum_{k=0}^{\infty}{z^k u^ke^{-ku}}du}\\
=\sum_{k=0}^{\infty} {z^{k+1}\int_0^{\infty}{u^k
e^{-(k+1)u}du}}=\sum_{k=1}^{\infty}{\frac{(k-1)!}{k^k}z^k}.
\end{eqnarray*}
We now use the theorem of Sibuya and Sperber \cite{sibuya}, which
gives the following necessary condition on a power series
$\sum_{k=0}^{\infty}{a_k z^k}$ with rational coefficients in order
for it to satisfy an algebraic differential equation:
\begin{equation}
\label{ncond} |a_k|_{p}\leq e^{Ck}\;\;\; for\; all\; k>0,
\end{equation}
where $|a|_{p}$ ($p$ prime) is the $p$-adic valuation defined by
writing $a=p^{i}\frac{m}{n}$ ($m,n$ not divisible by $p$) and
setting $|a|_{p}=p^{-i}$. Here we use this necessary condition
with $p=2$. We have $a_k=\frac{(k-1)!}{k^k}$. We now choose
$k=2^j$. The heighest power of $2$ dividing $(k-1)!$ is less than
$$\frac{k-1}{2}+\frac{k-1}{4}+\frac{k-1}{8}+...\leq k-1 = 2^j-1.$$
On the other hand the denominator of $a_k$ is $2^{j2^{j}}$. These
two facts imply that
$$|a_{2^j}|_2\geq 2^{(j-1)2^j+1},$$
which together with (\ref{ncond}) implies:
$$2^{(j-1)2^j+1}\leq |a_{2^j}|_2 \leq e^{C2^j},$$
or, taking logs
$$((j-1)2^j+1) log(2) \leq C2^j,$$
which is obviously false for $j$ sufficiently large. Hence the
necessary condition of Sibuya and Sperber is not satisfied for our
series, completing the proof.

\begin{figure}
\centering
    \includegraphics[height=8cm,width=9cm, angle=0]{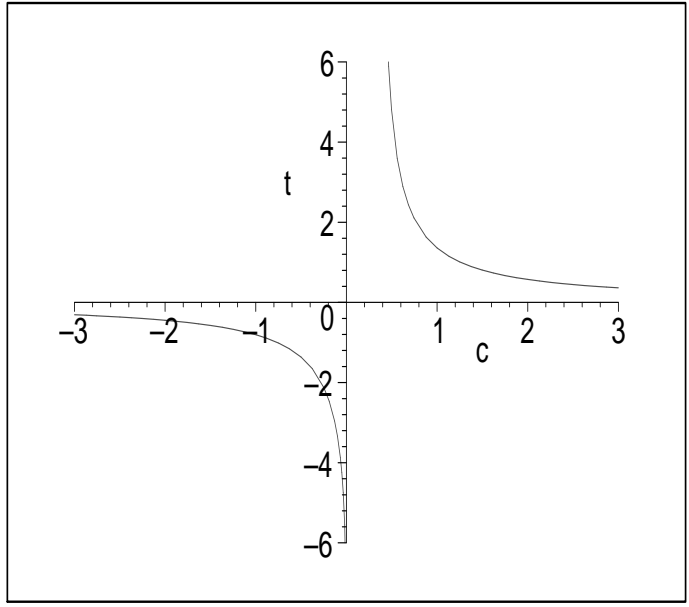}
    \caption{Graph of the function $H(c)$ defined by (\ref{defh})}
    \label{fone}
\end{figure}

\vspace{0.4cm} We note a few properties of the restriction of $H$
to the real line (which is defined outside the interval
$[0,\frac{1}{e}]$), obtained by elementary arguments. See figure
\ref{fone} for the graph of $H$, plotted using MAPLE.

\noindent (i) $H$ is decreasing on $(\frac{1}{e},\infty)$, with
$\lim_{c\rightarrow \frac{1}{e}+}{H(c)}=+\infty$,
$\lim_{c\rightarrow \infty}{H(c)}=0$.

\noindent (ii) $H$ is decreasing on $(-\infty,0)$, with
$\lim_{c\rightarrow -\infty}{H(c)}=0$, $\lim_{c\rightarrow
0-}{H(c)}=-\infty$.

\vspace{0.4cm}

We now define the function $F(s,c)$, where $s\in \Real$ is
considered a parameter and $c\in \Complex$ by:
\begin{equation}
\label{deff} F(s,c)=\int_0^s{\frac{du}{ce^{u}-u}}
\end{equation}

Since for $F(s,c)$ to be defined we need to ensure that the
denominator does not vanish for $u\in[0,s]$, it can be checked
that the domain of definition of $F$ is
$$D_F=\{ (s,c) \;|\; s\in \Real,\; c\in \Complex - J(s)\},$$
where
\begin{eqnarray*}
J(s)=[se^{-s},0]\;\;for\;s<0 \\
J(s)=[0,se^{-s}]\;\;for\;0\leq s<1\\
J(s)=[0,\frac{1}{e}]\;\;for\; s\geq 1.
\end{eqnarray*}

Using the fact that $H$ is {\it{not}} DA (lemma \ref{nda1}), we
now prove that the function $F(s,.)$ is {\it{not}} DA for generic
$s$. The argument is similar to the one in \cite{johnson} showing
that the function $z\rightarrow \int_{a}^{b}{u^ze^{-u}du}$ is not
DA for generic $a,b$ by using the fact that the Gamma function is
not DA.

\begin{lemma}
\label{nda} For generic $s\in \Real$, the function $F(s,.)$ is
{\it{not}} DA.
\end{lemma}

\noindent {\sc{proof:}} We apply lemma \ref{alternative} with
$p=1$, $\Field =\Complex$, $\Omega=D_F\subset \Real\times
\Complex$ (which is indeed open and connected). We have
$\Omega'=\Real$. We want to show that alternative (II) of lemma
\ref{alternative} holds. Let us assume by way of contradiction
that (I) holds, so that $F(s,.)$ is DA for all $s\in \Real$, and
in fact there is a common polynomial differential equation
satisfied by all $F(s,.)$:
\begin{equation}
\label{asd}
Q(c,F(s,c),D_cF(s,c),...,D_c^{(n)}F(s,c))=0,\;\;\;for\; all\;
(s,c)\in D_F.
\end{equation}
We now note that
$$\lim_{s\rightarrow\infty}{F(s,c)}=H(c)$$
uniformly on compact subsets of $D_H$. This implies that
derivatives of $D_c^{(k)}F(s,c)$ also converge to corresponding
derivatives of $H$, and thus from (\ref{asd}) it follows that
$$Q(c,H(c),H'(c),...,H^{(n)}(c))=0$$
for all $c\in D_H$, hence $H$ is DA, in contradiction with lemma
\ref{nda1}. This contradiction proves that alternative (II) holds,
as we wanted to show.

\vspace{0.4cm} We remark that we would guess that the
`exceptional' set in lemma \ref{nda} is $\{0 \}$, so that in fact
$F(s,.)$ is {\it{not}} DA for {\it{all}} $s\neq 0$, but we do not
know how to prove this.

\vspace{0.4cm}

We now restrict $c$ to be real nonzero number, and we define a new
function $G(t,c)$ by the relation:
\begin{equation}
\label{rel} G(F(s,c),c)=s\;\;\;for\; all\;(s,c)\in D_F.
\end{equation}
In other words we now look at $c$ as a parameter and define
$G(.,c)$ as the inverse function of $F(.,c)$. To see that $G$ is
well-defined and determine its domain of definition, we note the
following properties of $F$, which are elementary to verify.

\begin{lemma}
\label{props}

For each $c\leq \frac{1}{e}$ let $s^*(c)$ denote the solution (in
the case $0<c< \frac{1}{e}$, the {\it{smaller}} solution) of the
equation $se^{-s}=c$. We have:

\noindent (i) When $c<0$ the function $F(.,c)$ is decreasing on
$(s^*(c),\infty)$, \newline $\lim_{s\rightarrow
s^*(c)+}F(s,c)=+\infty$ and
$\lim_{s\rightarrow\infty}F(s,c)=H(c)<0$. \newline Hence $G(.,c)$
is defined on $(H(c),+\infty)$.

\noindent (ii) When $0<c\leq \frac{1}{e}$ the function $F(.,c)$ is
increasing on $(-\infty,s^*(c))$,\newline $\lim_{s\rightarrow
-\infty}F(s,c)=-\infty$ and $\lim_{s\rightarrow
s^*(c)-}F(s,c)=+\infty$. \newline Hence $G(.,c)$ is defined on
$(-\infty,\infty)$.

\noindent (iii) When $c> \frac{1}{e}$ the function $F(.,c)$ is
increasing on $(-\infty,\infty)$, \newline $\lim_{s\rightarrow
-\infty}F(s,c)=-\infty$ and $\lim_{s\rightarrow
\infty}F(s,c)=H(c)>0$. \newline Hence $G(.,c)$ is defined on
$(-\infty,H(c))$.

\vspace{0.4cm} Putting these facts together, we obtain that the
domain of definition $D_G$ of $G$ is:
$$D_G=D_G^{-}\cup D_G^{+},$$
where
$$D_G^{-}=\{ (t,c) \; |\; c<0,\; H(c)<t<\infty
\},$$
$$D_G^{+}=\{ (t,c) \; |\; 0<c\leq\frac{1}{e}\} \cup \{ (t,c)
\;|\; c>\frac{1}{e},\; -\infty<t<H(c)\}.$$ (in figure \ref{fone},
$D_G^{-}$ is the domain bounded by the left part of the graph of
$H$ and the $t$-axis, and $D_G^{+}$ is the domain bounded by the
right part of the graph and the $t$-axis).
\end{lemma}

We also note a fact that seems hard to prove directly, but which
follows indirectly from the results of the next section, as will
be pointed out.

\begin{lemma}
\label{cont} The function $G$ can be continued as a real-analytic
function to the open connected domain
\begin{equation}
\label{dd} D=D_G^{-}\cup D_G^{+}\cup \{ (t,0) \; |\; t\in \Real \}
\end{equation}
by setting
\begin{equation}
\label{ext} G(t,0)=0 \;\;\;for\; all\; t\in \Real
\end{equation}
(in figure \ref{fone}, $D$ is the domain bounded between the two
parts of the curve representing the graph of $H$).
\end{lemma}

\begin{lemma}
\label{invda} For generic $t\in \Real$, the function $G(t,.)$ is
{\it{not}} DA.
\end{lemma}

\noindent {\sc{proof:}} Assume by way of contradiction that the
conclusion of the lemma is not true, so that $G(t,.)$ is DA for a
set of values of $t$ which is either of positive measure or not of
first category. Then, applying lemma \ref{alternative}, with
$\Omega=D$ (recall (\ref{dd})), we conclude that $G(t,.)$ is DA
for {\it{all}} $t\in \Real$. We note also that $G(.,c)$ is DA for
all value of $c$, since for $c\neq 0$ it is the inverse of
$F(.,c)$ which is DA, and the inverse of a DA function is DA,
while for $c=0$ it is identically $0$ by (\ref{ext}). Thus $G$ is
DA with respect to both variables on $\Omega$. Now fixing any
$\overline{s}\in \Real$, and using the defining relation
(\ref{rel}) we have
$$G(F(\overline{s},c),c)-\overline{s}=0$$
for all $c$. We now use lemma \ref{implicit}, with
$R(w,z)=G(w,z)-\overline{s}$, to conclude that $F(\overline{s},.)$
is DA. Since $\overline{s}\in \Real$ was aribtrary, we have that
$F(s,.)$ is DA for all $s\in \Real$, but this contradicts the
result of lemma \ref{nda}. This contradiction concludes the proof.

\section{Construction of a differential equation with non-DA
dependence on initial conditions}

\label{const}

\label{example}

\begin{theorem}
\label{difeq} Consider the initial value problem:
\begin{eqnarray}
\label{eq1}
y_{1}'=y_{2}-y_{1},\;\; y_{1}(0)=r_1\\
\label{eq2} y_{2}'=y_{2}(y_{2}-y_{1}),\;\; y_{2}(0)=r_2
\end{eqnarray}
Let the solutions be denoted by $y_1(r_1,r_2;x)$,
$y_2(r_1,r_2;x)$, defined on the maximal interval $I(r_1,r_2)$.
Then for generic $x\in \Real$, the functions $y_1(0,r_2;x)$,
$y_2(0,r_2;x)$ are not DA with respect to $r_2$.

\end{theorem}

To prove theorem \ref{difeq}, we assume for the moment that
$r_2\neq 0$, which implies that $y_2(x)\neq 0$ for all $x$ (since
by (\ref{eq2}) and by the standard uniqueness theorem for ODE's,
if $y_2$ vanishes at one point, it vanishes identically), and note
that from (\ref{eq1}) and (\ref{eq2}) we have
$$\frac{y_2'}{y_2}=y_1'$$
which implies
$$y_2(x)=y_2(0)e^{(y_1(x)-y_1(0))},$$
and in the case $y_1(0)=r_1=0, \;y_2(0)=r_2$:
\begin{equation}
\label{relate} y_2(x)=r_2e^{y_1(x)}.
\end{equation}
Substituting (\ref{relate}) back into (\ref{eq1}), we have
$$y_1'(x)=r_2e^{y_1(x)}-y_1(x),$$
or,
$$\frac{y_1'(x)}{r_2e^{y_1(x)}-y_1(x)}=1.$$
Thus, integrating, and recalling the function $F(s,c)$ defined by
(\ref{deff}), we have:
$$F(y_1(0,r_2;x),r_2)=x$$
for all $x\in I(0,r_2)$. In other words, recalling the definition
(\ref{rel}) of $G(t,s)$, we have:
\begin{equation}
\label{solution} y_1(0,r_2;x)=G(x,r_2)\;\;\; for\; all\; r_2\neq
0,\; x\in I(0,r_2)
 \end{equation}

We note in passing that, since $y_1(0,r_2;x)$ is real-analytic in
its domain of definition, and is defined also for $r_2=0$, with
$y_1(0,0;x)=0$, (\ref{solution}) implies that $G$ can be extended
to $D$ as lemma \ref{cont} claimed.

In lemma \ref{invda} we showed that $G(t,c)$ is {\it{not}} DA with
respect to $c$ for generic $t$. From (\ref{solution}) we then get,
for generic $x$, that $y_1(0,r_2;x)$ is {\it{not DA}} with resepct
to $r_2$. Using this fact and (\ref{relate}), together with the
fact that compositions and products of DA function are DA, we get
the same conclusion for $y_2$, concluding the proof of theorem
\ref{difeq}.

\section{proof of the main theorem}

\label{proof}

We first set some notation. Define $\Sigma \subset {\bf{S}}(m,d)
\times \Real \times \Real^m$ by
$$\Sigma= \{ (P,x,r)\;|\; P\in {\bf{S}}(m,d),\; r\in
\Real^m,\; x\in I(P,r)\}.$$ For each $1\leq j\leq m$ let
$$\Sigma'_j=\{ (P,x,r_1,...,r_{j-1},r_{j+1},...,r_{m}) \;|\; \exists r_j\;
\; such\;that\; (P,x,r)\in \Sigma \}.$$

The proof of theorem \ref{main}, will be based on the following

\begin{lemma}
\label{generic} Assume $m\geq 2$, $d\geq 2$. Fix $1\leq j \leq m$.
\newline Then for generic $(P,r_1,...,r_{j-1},r_{j+1},...r_m,x)\in
\Sigma'_j$, the functions $f_{ij}$ ($1\leq i\leq m$) defined by
(\ref{fi}) are {\it{not}} DA.
\end{lemma}

\noindent {\sc{proof:}} We note first that we may, without loss of
generality and for notational convenience, assume that $j=m$.
Define $F_k:\Sigma\rightarrow \Real$ ($1\leq k \leq m$) by
$$F_k(P,x,r_1,...,r_m)=y_k(r_1,...,r_m;x).$$
These are real analytic-function on $\Sigma$, and we may apply
lemma \ref{alternative}, with $\Field=\Real$, $\Omega=\Sigma$, $u$
identified with $(P,x,r_1,...,r_{m-1})$, $v$ identified with
$r_m$, and $F$ being any of the $F_k$ ($1\leq k \leq m$). We would
like to rule out alternative (I). We now fix some $1\leq i \leq
m-1$, and note that among the systems $P\in {\bf{S}}(m,d)$ is the
system
$$y_i'=y_m-y_i$$
$$y_m'=y_m(y_m-y_i)$$
$$y_k'=0 \;\;for\;k\neq i,m.$$
By theorem \ref{difeq}, for this system, for $r_i=0$ and for
arbitrary values of $r_k$ ($k\neq i,m$), the dependence of
$y_i(r_1,...,r_m;x)$ and of $y_m(r_1,...,r_m;x)$ on $r_m$ is not
DA. This rules out alternative (I) of lemma \ref{alternative} both
for $F_i$ and for $F_m$, which implies that alternative (II) holds
for both. Since $1\leq i\leq m-1$ is arbitrary, this  means that
for each $1\leq k \leq m$ there is a generic subset $\Sigma^g_k$
of $\Sigma'_m$ for which the function $F_k$ is {\it{not}} DA with
respect to $r_m$. Thus, setting
$$\Sigma^g=\cap_{1\leq k \leq m} \Sigma^g_k$$
we get a generic set, and for each $(P,x,r_1,...r_{m-1})\in
\Sigma^g$ we have the conclusion of lemma \ref{generic}.

\vspace{0.4cm} Finally, to derive theorem \ref{main} from lemma
\ref{generic}, we recall two classical results (see
\cite{oxtoby}). Let $C\subset \Real^N\times \Real^M$, and define
for each $u\in \Real^N$,
$$C(u)=\{ v\in \Real^M \;|\; (u,v)\in C\}.$$
We have:

\noindent (i) If the set $C$ is of measure $0$ then $C(u)$ is of
measure $0$ in $\Real^M$ for all $u\in \Real^N$ except a set of
measure $0$. This is (a special case of) Fubini's theorem.

\noindent (ii) If the set $C$ is of first category then $C(u)$ is
of first category in $\Real^M$ for all $u\in \Real^N$ except a set
of first category. This is the Banach-Kuratowski theorem.

\vspace{0.4cm}

\noindent {\sc{proof of theorem \ref{main}}:} We will prove the
conclusion of theorem \ref{main} for all $P\in {\bf{S}}(m,d)$
outside a set of measure $0$, making use of Fubini's theorem. To
prove the same conclusion for $P$ outside a set of first category,
one only needs to replace Fubini's theorem with the
Banach-Kuratowski theorem.

For $1\leq i,j \leq m$ define $C_{ij}$ to be the set of
$(P,x,r_1,...,r_{j-1},r_{j+1},...,r_{m})\in \Sigma_j'$  for which
the function $f_{ij}$ defined by (\ref{fi}) is DA. By lemma
\ref{generic}, each of the sets $C_{ij}$ is of measure $0$ and of
first category. Hence by Fubini's theorem, for all $P\in
{\bf{S}}(m,d)$ except a set of measure $0$, the set
$$C_{ij}(P)=\{(x,r_1,...,r_{j-1},r_{j+1},...,r_{m})\;|\;
(P,x,r_1,...,r_{j-1},r_{j+1},...,r_{m})\in C_{ij}\}$$ has measure
$0$. Below we show that the sets $C_{ij}(P)$ are also of first
category for any $P$ for which they are of measure $0$. Since this
conclusion is valid for any $1\leq i,j\leq m$, we have that for
all $P\in {\bf{S}}(m,d)$ except a set of measure $0$, {\it{all}}
the sets $C_{ij}(P)$ have measure $0$ and are of first category.
This is precisely the content of theorem \ref{main}.

It is left to show that whenever $C_{ij}(P)$ is of measure $0$, it
is also of first category. Without loss of generality, and for
notational convenience, we assume $j=m$. We apply lemma
\ref{alternative} with $\Field=\Real$, $\Omega=\Lambda(P)$ (recall
(\ref{lambda})), $u$ identified with $(x,r_1,...,r_{m-1})$, $v$
identified with $r_m$, and $F=f_{im}$. Since we assume that
$C_{im}(P)$ is of measure $0$, so that $F$ is {\it{not}} DA with
respect to $r_m$ for almost all $(x,r_1,...,r_{m-1})\in
\Lambda_{m}(P)$, alternative (I) of lemma \ref{alternative}
certainly does not hold, so (II) holds which means that
$C_{im}(P)$ is indeed also of first category, concluding the
proof.


\begin{thebibliography}{9}

\bibitem{ritt}  E. Gourin \& J. F. Ritt, An assemblage-theoretic proof of
 the existence of transcendentally transcendental functions,
 Bull. Amer. Math. Soc. {\bf{33}} (1927), 182-184.


\bibitem{johnson} J. Johnson, G.M. Reinhart \& L.A. Rubel, Some
counterexamples to seperation of variables, J. Diff. Eq.
{\bf{121}} (1995), 42-66.

\bibitem{oxtoby} J. Oxtoby, `Measure and Category',
Springer-Verlag (New-York), 1971.

\bibitem{pour} M.B. Pour-El \& J.I. Richards, `Computability in Analysis and
Physics', Springer-Verlag (Berlin), 1989.


\bibitem{rubel4} L.A. Rubel, {\it{An elimination theorem for systems of
 algebraic differential equations}}, Houston J. Math. {\bf{8}} (1982), 289-295.

\bibitem{rubelf} L.A. Rubel, {\it{Some research problems about
algebraic differential equations}}, Trans. Am. Math. Soc.
{\bf{280}} (1983), 43-52.

\bibitem{rubel2} L.A. Rubel \& M.F. Singer, {\it{A differentially algebraic elimination
theorem with application to analog computability in the calculus
of variations}}, Proc. Am. Math. Soc. {\bf{94}} (1985), 653-658.

\bibitem{rubel3} L.A. Rubel, {\it{A survey of transcendentally transcendental
functions}}, Amer. Math. Monthly {\bf{96}} (1989), 777-788.

\bibitem{rubel} L.A. Rubel, {\it{Some research problems about
algebraic differential equations II}}, Illinois J. Math. {\bf{36}}
(1992), 659-680.

\bibitem{sibuya} Sibuya, Yasutaka \& Sperber, Steven, Arithmetic properties of power series solutions of algebraic differential
equations, Ann. of Math. {\bf{113}} (1981), 111--157.

\end{thebibliography}
\end{document}